\newlength{\extralineskip}

  \documentstyle[12pt]{article}
\input epsf
\newcommand{\be}{\begin{equation}}
\newcommand{\ee}{\end{equation}}
\newcommand{\bae}{\begin{eqnarray}}
\newcommand{\eae}{\end{eqnarray}}
\def\men#1{~(\ref{#1})}
\def\mens#1{~\ref{#1}}
\newtheorem{prop}{Proposition}
\newtheorem{lemma}{Lemma}
\def\fgr#1#2{\hbox{\raise #1 ex\hbox{\epsfbox{#2}}}}
\newlength{\lpr}
\newlength{\wpr}
\def\prend{$\Box$}
\def\bb{\bibitem}

\def\ie{\hbox{\it i.e.\,}}
\def\etc{\hbox{\it etc.\,}}
\def\eg{\hbox{\it e.g.\,}}
\def\nn{\nonumber}
\def\ff{\nn \\}
\def\smash{{\A \rtimes \U}}
\def\GL{GL_{q}(N)}
\def\A{{\cal A}}

\def\D{{\cal D}}

\def\I{{\cal I}}

\def\P{{\cal P}}
\def\R{{\cal R}}

\def\U{{\cal U}}

\def\X{{\cal X}}

\def\xisi{\xi_{i_{1}} \ldots \xi_{i_{k}} \uot \sigma_{i_{k}} \ldots
\sigma_{i_{1}}}
\def\xs{\xi_{i_{1}} \ldots \xi_{i_{k}}}
\def\sx{\sigma_{i_{k}} \ldots \sigma_{i_{1}}}

\def\Ldu{\delta^{\scr L}_{\cal U}}
\def\Rdu{\delta^{\scr R}_{\cal U}}
\def\da{\delta_{\cal A}}
\def\dX{\delta_{\cal X}}
\def\Rda{\delta^{\scr R}_{ \cal A}}
\def\Rdx{\delta^{\scr R}_{ \cal X}}
\def\Lda{\delta^{\scr L}_{ \cal A}}
\def\Ldx{\delta^{\scr L}_{ \cal X}}
\def\p{\partial}
\def\E{{\cal E}}
\def\bE{\bar{\E}}

\def\hR{\hat{R}}
\def\delA{\Delta_{\A}}

\def\Udel{_{\U}\Delta}
\def\Sfi{S^{2}(f^{i}}
\def\Sfk{S^{2}(f^{k}}
\def\Sfj{S^{2}(f^{j}}
\def\Sioo{S^{-1}({a^{(\bar{1})}}_{(1)})} 
\def\Saot{S^{-1}({a^{(\bar{1})}}_{(2)})} 
\def\Stta{S^{-1}({a^{(2)}}_{(2)})} 
\def\Sia{S^{-1}(a^{(\bar{1})}} 
\def\id{\mathop{\rm id}}

\def\uA{1_{\A}}

\def\uU{1_{\U}}
\def\ot{\otimes}
\def\uot{\underline{\otimes}}
\newcommand{\rtimes}{\mbox{$\times\!\rule{0.3pt}{1.1ex}\,$}}
\def\tr{\triangleright}
\def\ra{\rightarrow}

\def\fe{& = &}

\def\detq{\mathop{\rm det_{q}}}

\def\vaca{\Omega_{\A}}
\def\lva{\langle \Omega_{\A} |}
\def\rva{| \Omega_{\A} \rangle}
\def\vacu{\Omega_{\U}}
\def\ulv{\langle \Omega_{\U} |}
\def\urv{| \Omega_{\U} \rangle}
\def\la{\langle}  
\def\ar{\rangle}
\def\Rar{\ar^{\scr R}}
\def\Lar{\ar^{\scr L}}
\def\trRar{\Rar_{tr}}
\def\turRar{\Rar_{\underline{tr}}}

\def\inprod#1#2{\left\langle #1, #2\right\rangle}
\def\ip#1#2{\left\langle #1, #2\right\rangle}
\def\olo#1{\langle #1 \rangle}
\newcommand{\bolo}[1]{\langle\!\langle {#1} \rangle \! \rangle}

\def\Lbolo#1{\langle\!\langle #1  \rangle  \! \rangle^{\scr  L}}
\def\Rolo#1{\langle #1 \Rar}
\def\Lolo#1{\langle #1 \Lar}
\def\ftR#1{\widehat{#1}^{\scr R}}

\def\scr{\scriptscriptstyle}

\def\sigmai{\sigma_{i}}

\def\p{\partial}

\def\udel{\underline{\Delta}}

\begin{document}
\newcommand{\norm}[1]{{\protect\normalsize{#1}}}
\newcommand{\LAP}
{{\small E}\norm{N}{\large S}{\Large L}{\large A}\norm{P}{\small P}}
\newcommand{\sLAP}{{\scriptsize E}{\footnotesize{N}}{\small S}{\norm L}$
${\small A}{\footnotesize{P}}{\scriptsize P}}
\begin{titlepage}
\begin{minipage}{5.2cm}
\begin{center}
{\bf Groupe d'Annecy\\
\ \\
Laboratoire d'Annecy-le-Vieux de Physique des Particules}
\end{center}
\end{minipage}
\hfill
\hfill
\begin{minipage}{4.2cm}
\begin{center}
{\bf Groupe de Lyon\\
\ \\
Ecole Normale Sup\'erieure de Lyon}
\end{center}
\end{minipage}
\centerline{\rule{12cm}{.42mm}}

\vskip 2cm
\begin{center}

  {\LARGE {\bf Remarks On Quantum Integration}}\\[1cm]

  \vskip 2cm

  {\large Chryssomalis Chryssomalakos$^{\ast}$}

  \vskip 1cm

  {\em
    Laboratoire de Physique Th\'eorique }\LAP\footnote{URA 14-36
    du CNRS, associ\'ee \`a l'Ecole Normale Sup\'erieure de Lyon et \`a
    l'Universit\'e de Savoie.

    \indent
    $^\ast$chryss@lapp.in2p3.fr
    }\\
    {\em Chemin de Bellevue BP 110, F-74941 Annecy-le-Vieux Cedex,
      France.}

\end{center}
\vskip 2cm
\begin{abstract}
We give a general integration prescription for finite dimensional
braided Hopf algebras, deriving the N-dimensional quantum
superplane integral as an example. The transformation properties of
the integral on the quantum plane are found. We also discuss 
integration on quantum group modules that lack a Hopf structure.

\end{abstract}

\vskip 5cm

\hfill \LAP -A-562/95

\hfill December 1995

\end{titlepage}
\newpage
\pagestyle{plain}
\section{Introduction}
The emergence of Hopf algebras, during the last decade, as a
promising framework within which new physical symmetries can be
accomodated, has prompted an interest in the theory and techniques
of integration on them. Similar remarks hold for braided Hopf
algebras which, more recently, have provided a still further
extension of the classical concept of a group, marrying
quantization with nontrivial statistics. 

Integrals on (finite dimensional) Hopf algebras have been studied
extensively, see for example~\cite{LarSwee,Rad} and references
therein. For the braided case see the treatment in~\cite{Lyu1,Lyu2}
for the basics of the theory - some examples appear
in~\cite{CZ,MajKem}. From the point of view of a physicist who is
interested in the basic elements of the theory and in applications,
the situation presernts certain problems. The results are generally
scattered and when they do become available, the disparity of the
methods employed in them prohibits the formation of a clear image
of the minimum background required to explore the field. When it
comes to applications, results are extracted (often ingeniously)
from particular properties of individual examples - no general
integration prescription seems to be available for the quantum
space wanderers. The closest one can come to such a prescription in
the literature is perhaps the trace formula of Radford and
Larson~\cite{LarRad1} (and a braided analogue of it in~\cite{Majid6})
which, however, occasionally returns trivial (\ie identically zero)
results (there doesn't seem to exist a description of when exactly
it fails either). Our purpose in this paper is thus basically
twofold. On the one hand, we aim at providing simple, self
contained proofs of basic results, using only the Hopf algebra
axioms we assume the reader to be familiar with, attempting in this
way a demonstration of what can be accomplished with a rather
minimal set of tools. On the other hand, addressing the problem of
the missing integration prescription, we give an explicit formula
for the integral on any finite dimensional braided Hopf algebra
(FDBHA) and show that it is always nontrivial, commenting along the
way on the conditions under which the trace formula fails.

The paper is structured as follows: in section\mens{HA} we present
the notation we use and collect some basic formulas we need in
subsequent proofs. Section\mens{IoHA} starts with background
information on integrals on (nonbraided) Hopf algebras. We then
give a modified trace formula for the integral and prove its
nontriviality. We also introduce a ``vacuum expectation value''
approach to integration, discuss properties of right Fourier
transforms and prove a number of useful formulas.
Section\mens{IoBHA} supplies the braided version of the modified
trace formula and the vacuum projectors and discusses, as an
example, the integral on the N-dimensional quantum superplane. Also
included are some comments on the transformation properties of the
integral on the N-dimensional quantum plane. The last section
provides an integration prescription (with some modest assumptions)
for quantum group modules that lack a braided Hopf structure.
\section{Hopf Algebras} \label{HA}
The language used in the following is predominantly that of Hopf
algebras - we refer the reader to~\cite{Abe,Majid1,SweeHA} for an
introduction to the subject. Concerning the notation,
we denote by $\Delta$, $\epsilon$, $S$ the {\em coproduct}, {\em
 counit} and {\em antipode} respectively and by $\delA$, $\Udel$
the {\em right} $\A$ and {\em left} $\U$-{\em coactions}
respectively. Sweedler-like conventions are employed throughout -
thus $\Delta{a} = a_{(1)} \ot a_{(2)}$, $(\Delta \ot \id) \circ
\Delta(a) = a_{(1)} \ot a_{(2)} \ot a_{(3)}$ \etc . Also, $\delA
(x) = x^{(1)} \ot x^{(2')}$, $\Udel(a) = a^{(\bar{1})} \ot a^{(2)}$
and, for example, $(\id \ot \Delta) \circ \delA (x) = x^{(1)} \ot
{x^{(2')}}_{(1)} \ot {x^{(2')}}_{(2)} =  x^{(1)} \ot
x^{(2')} \ot x^{(3')}$ and so on. By $\A$ we will generally denote a
function type Hopf algebra (its elements will be denoted by $a$,
$b$, \etc) - $\U$ will stand for its dual Hopf
algebra (universal enveloping algebra type) with elements $x$, $y$
\etc . The duality is via a nondegenerate {\em inner product}
$\ip{\cdot}{\cdot}$ that relates the algebra structure in $\A$ with the
coalgebra stucture in $\U$ and vice-versa. The {\em universal}
$R$-{\em matrix} is denoted by $\R = \R^{(1)} \ot \R^{(2)}$; $\R'$
stands for $\tau (\R)$ with $\tau(a \ot b) = b \ot a$.
$\R$ satisfies
\be
\Delta ' (x) = {\cal R} \Delta (x) {\cal R}^{-1}  \qquad \qquad
\forall x \in \U 
\label{RcopRi}
\ee
as well as
\bae
(\Delta \otimes \id ){\cal R} \fe {\cal R}^{13} {\cal R}^{23} 
\label{Rprop1} \\
( \id \otimes \Delta) {\cal R} \fe {\cal R}^{13} {\cal R}^{12} 
\label{Rprop2}
\eae
(with ${\cal R}^{13} \equiv {\cal R}^{(1)} \otimes \uU \otimes 
{\cal R}^{(2)}$ {\em etc.}) 

Given a pair
of dual Hopf algebras, one can construct their {\em semidirect
product} $\smash$ with $\A$, $\U$ trivially embedded in it
and cross relations
\be
xa = a_{(1)} \inprod{x_{(1)}}{a_{(2)}} x_{(2)}  \, ,
\qquad
ax = x_{(2)} \inprod{x_{(1)}}{S^{-1}(a_{(2)})} a_{(1)} \, .
\label{AUcr}
\ee
The above commutation relations guarrantee that~\cite{VD2}
\be
\A \ot \U \ni a \ot x \neq 0 \Rightarrow xa \neq 0 \, 
\label{VD1}
\ee
with $xa \in \smash$. The {\em action} $\tr$ of $\U$ on $\A$ is
 given by $x \tr a =
a_{(1)} \ip{x}{a_{(2)}}$. The same symbol will denote the ({\em
adjoint}) {\em action} of
$\U$ on $\U$: $x \tr y = x_{(1)} y S(x_{(2)}) = y^{(1)}
\ip{x}{y^{(2')}}$. The canonical element in $\U \ot \A$ is written
like $ C = e_{i} \ot f^{i}$ with $\{e_{i}\}$, $\{f^{i}\}$ dual
(in the sense that $\ip{e_{i}}{f^{j}} = \delta_{i}^{j}$) 
linear bases in $\U$, $\A$ respectively. It holds
\bae
(\Delta \ot \id ) C_{12} \fe C_{13} C_{23}  \ff
(\id \ot \Delta ) C_{12} \fe C_{12} C_{13}  \ff 
(\epsilon \ot \id) C \fe 1 \ff
(\id \ot \epsilon) C \fe 1 \ff
(S \ot \id) C \fe C^{-1} \nn \\
(\id \ot S) C \fe C^{-1} \, , 
\label{Cprop1}
\eae
as well as 
\bae
\delA(a) \equiv \Delta(a) \fe C(a \ot 1) C^{-1} \ff
\delA(x) \fe C(x \ot 1) C^{-1} \ff
\Udel(a) \fe C^{-1} (1 \ot a) C \ff
\Udel (x) \equiv \Delta(x) \fe C^{-1} (1 \ot x) C  \, . 
\label{Cprop2}
\eae
Either of\men{AUcr} encodes the information about the
inner product $\inprod{x}{a}$. To make this precise, we
introduce $\U$ and $\A$-{\em right vacua}, denoted by $\urv$ and
$\rva$ respectively, which satisfy~\cite{CZ}
\bae
x \urv \fe \epsilon(x) \urv \ff
a \rva \fe \epsilon(a) \rva \nn .
\eae
{\em Left vacua} $\ulv$, $\lva$ are defined analogously.
In terms of these, the inner product $\inprod{x}{a}$ can be given
as the ``expectation value''
\bae
\lva xa \urv \fe \lva a_{(1)} \inprod{x_{(1)}}{a_{(2)}} x_{(2)}
\urv \ff
 \fe \inprod{x}{a}
\eae
if we normalize the vacua so that $\langle \vaca | \vacu \rangle =
\langle \vacu | \vaca \rangle = 1$. Similarly, the adjoint action
of $\U$ on $\A$ can be written as
\be
xa \urv = x \tr a \urv \nn  .
\ee
\section{Integration on Hopf Algebras} \label{IoHA}
\subsection{Background} \label{B}
We list here known results about invariant integrals on Hopf
algebras that we will use later - more details can be found
in~\cite{SweeHA,Abe}. Some of the proofs are also
supplied in order to familiarize the reader with the usage of the
tools presented in section\mens{HA}. To prevent potential
divergence problems from distracting the formulation of concepts,
we deal throughout with finite dimensional Hopf algebras - some of
the results though retain their validity in the infinite
dimensional case as well.

We start with two dually paired Hopf algebras $\U$, $\A$ with
(finite) dual bases $\{e_{i}\}, \; \{f^{j}\}$ respectively. We define a
{\em right (invariant) integral} in $\A$ as a map $\la \cdot
\Rar$: $\A \ra k$ with the property
\be
\la a_{(1)} \Rar a_{(2)} = \la a \Rar \uA \label{rint}
\ee
for all $a$ in $\A$. We call $\la \cdot \Rar$ {\em trivial} if all
$\la f^{i} \Rar$ are zero. {\em Left (invariant) integrals} are
similarly defined via
\be
a_{(1)} \la a_{(2)} \Lar = \uA \la a \Lar . \label{lint}
\ee
As we shall soon see, when $\la \uA \Rar \neq 0$ (or $\la \uA \Lar
\neq 0$), left and right integrals are proportional and
can therefore be normalized so that they coincide.
One can now introduce the element $\Rdu \in \U$ (the {\em right
delta function} in $\U$) which 
implements the right integral in $\A$ via
\be
\inprod{\Rdu}{a} = \la a \Rar 
\ee
(so that $\epsilon(\Rdu)= \la \uA \Rar $). This allows us to write
\be
\Rdu = \la f^{i} \Rar e_{i} 
\ee
(in the mathematics literature $\Rdu$ is often called a {\em
right integral} in $\U$ - we will not use this terminology here).
For arbitrary $a$ in $\A$ we have
\bae
\inprod{\Rdu x}{a} \fe \inprod{\Rdu}{a_{(1)}} \inprod{x}{a_{(2)}} \ff
 & = & \la a_{(1)} \Rar \inprod{x}{a_{(2)}} \ff
 & = & \la a \Rar \inprod{x}{\uA} \ff
 & = & \epsilon(x) \la a \Rar \ff
 & = & \inprod{\Rdu \epsilon(x) }{a} \nn
\eae
therefore
\be
\Rdu x = \Rdu \epsilon(x)  \; \; \; \; \; \; \forall x \in \U.
\label{xdU}
\ee
Taking antipodes in the above equation we find ($\epsilon(x) =
\epsilon(S(x))$)
\be
S(x) S(\Rdu) = \epsilon(S(x)) S(\Rdu) \nn
\ee
which, for invertible $S$ gives
\be
x S(\Rdu) = \epsilon(x) S(\Rdu) \; \; \; \forall x \in \U , \nn
\ee
in other words, $S(\Rdu)$ implements a left integral and
we can therefore take $S(\Rdu) = \Ldu$. Since $\epsilon(S(x)) =
\epsilon(x)$, we find that $\la \uA \Rar = \la \uA \Lar \equiv \la
\uA \ar$ ($\olo{\cdot}$ denotes a biinvariant integral). 
Consider now the product $\Rdu S(\Rdu)$. We have
\[
\Rdu S(\Rdu) = \Rdu \epsilon(S(\Rdu)) = \Rdu \epsilon(\Rdu) = \Rdu
\la \uA \ar
\]
and also
\[
\Rdu S(\Rdu) = \epsilon(\Rdu) S(\Rdu) = \la \uA \ar S(\Rdu) 
\]
therefore, in the unimodular case where $\la \uA \ar \neq 0$ (so
that we can normalize $\olo{\uA}=1$), it holds 
\bae
\Rdu \fe S(\Rdu) \label{duSdu} \\
x \Rdu \fe \epsilon(x) \Rdu \; \; \forall x \in \U \label{dux} \\
\la a \ar \fe \la S(a) \ar \; \; \forall a \in \A  \label{intaSa}
\, .
\eae
Concerning uniqueness, assume that a second right
integral $\la \cdot {\Rar}^{'}$ exists and let ${\Rdu}^{'}$ be
 the element
of $\U$ that implements it. We then get $\Rdu S({\Rdu}^{'}) =
\la \uA \ar S({\Rdu}^{'})$ and $\Rdu S({\Rdu}^{'}) = \Rdu \la \uA
\ar^{'} $. For $\la \uA \ar =1$ we conclude that $\la \uA \ar^{'} =
1 $ implies $\Rdu = S({\Rdu}^{'}) = {\Rdu}^{'}$ while $\la \uA \ar^{'}
= 0$ implies that $\la \cdot \ar^{'}$ is trivial.

\noindent Radford and Larson~\cite{LarRad1} have shown that
\be
\la a \trRar \equiv \ip{S^{2}(e_{i})}{f^{i}a} \label{RadL}
\ee
defines a right integral on $\A$ which though, in some
cases (as~\cite{LarRad1} warns), is trivial. It is this shortcoming
of\men{RadL} that (among other things) motivated our formula for
the integral of the next section, which is shown to be 
nontrivial for any FDHA. One can derive from\men{RadL}
that
$\inprod{S^{-2}(e_{i})}{ f^{i}a}$,
$\inprod{S^{-2}(e_{i})}{ af^{i}}$ define left and right
(again, in some cases trivial) integrals respectively. For $\sigma
\equiv \ip{S^{2}(e_{i})}{f^{i}} \neq 0$ and $\sigma' \equiv
\ip{S^{-2}(e_{i})}{f^{i}}$ we conclude
\be
\inprod{S^{-2}(e_{i})}{f^{i} a} = \inprod{S^{-2}(e_{i})}{a f^{i}} =
\frac{\sigma'}{\sigma} \la a \trRar . \label{lrint2} 
\ee
We summarize the main points: when $\la \uA \ar = 1$, $\la \cdot
\ar$ is biinvariant, unique and\men{intaSa}
holds. \men {RadL} defines a (sometimes trivial) right 
integral and when $\sigma \neq 0$,\men{lrint2} holds. Considerably
more is known about integrals on Hopf algebras - in the interest of
self containment we have only mentioned above what we can prove
here. 
\subsection{A modified trace formula}
\label{altform}
We want to present now a modified version of the trace
formula\men{RadL} that overcomes the limitations mentioned above.
It is given by 
\be
\la a \Rar \Rda = \ip{e_{j} S^{-2}(e_{i})}{a f^{i}} f^{j} 
\label{intdef3}
\ee
(for the remainder of this paper, $\la a \Rar \Rda$ is 
defined by the rhs
of\men{intdef3}). Notice that\men{intdef3} also defines (for some
nonzero $\la a \Rar$ which we can normalize to 1) the right delta
function in $\A$. By pairing both sides of\men{intdef3} with $x$ in
$\U$ we conclude
\be
\la x \Rar \la a \Rar = \ip{x S^{-2}(e_{i})}{a f^{i}} .
\label{intdefxLaR}
\ee
The proof of invariance is quite analogous to that of\men{RadL}.
What is interesting is the following
\begin{lemma}
The integral $\la \cdot \Rar$ defined by\men{intdef3} is nontrivial
for any FDHA $\A$.
\end{lemma}
{\em Proof:} Set $\Theta^{k}_{l} \equiv \ip{e_{l} S^{-2}(e_{i})}{f^{k}
f^{i}}$ and compute
\bae
\smash \ni S^{-2}(e_{i}) f^{i} \fe f^{i}_{(1)}
\ip{S^{-2}(e_{i_{(1)}})}{f^{i}_{(2)}} S^{-2}(e_{i_{(2)}}) \ff
 \fe f^{i}_{(1)} f^{j}_{(1)} \ip{S^{-2}(e_{i})}{f^{i}_{(2)}
f^{j}_{(2)}} S^{-2}(e_{j}) \ff
 \fe S^{-2}(f^{i}) f^{j} \ip{e_{i} S^{-2}(e_{k})}{f^{k} f^{l}}
S^{-2}(e_{j}) S^{-2}(e_{l}) \ff
 \fe \Theta^{l}_{i} S^{-2}(f^{i}) f^{j} S^{-2}(e_{j}) S^{-2}(e_{l}).
\eae
We now employ\men{VD1} to conclude that not all $\Theta^{l}_{i}$
are zero (since $S^{-2}(e_{i}) \otimes f^{i} \neq 0$ in $\U \otimes
\A$). Alternatively, one can compute directly the
integral~\cite{VD1}
\bae
\Rolo{S^{-1}(\Rda)} \fe \Rolo{e_{i}} \Rolo{S^{-1}(f^{i})} \ff
 \fe \ip{e_{i} e_{j}}{S^{-1}(f^{i}) S^{-2}(f^{j})} \ff
 \fe \ip{e_{i}}{S^{-1}(f^{i}_{(1)}) S^{-2}(f^{i}_{(2)})} \ff
 \fe \ip{e_{i}}{\epsilon(f^{i}) \uA} \ff
 \fe \ip{\uU}{\uA} \ff
 \fe 1 \, , \label{intdeltaR}
\eae
which shows nontriviality in both $\U$ and $\A$ (of course,
$\Rolo{e_{i}} \Rolo{S^{-1}(f^{i})} = 1$ is a stronger statement).
\hfill \prend

\noindent Defining $\Theta = \Theta^{(1)} \otimes \Theta^{(2)} =
\Theta^{l}_{k} f^{k} \otimes e_{l}$ we get
\be
\Rolo{a} \Rda = \Theta^{(1)} \ip{\Theta^{(2)}}{a} \, \qquad
\Rolo{x} \Rdu = \ip{x}{\Theta^{(1)}} \Theta^{(2)} \, .
\label{axTheta}
\ee

\noindent Related to the above proof is
\begin{lemma}
For $\A$ a FDHA, it holds ($a$ in $\A$)
\be
\Rolo{af^{i}} = 0 \, \, \, \, \forall i \Rightarrow a=0
\label{intnondeg}
\ee
\end{lemma}

{\em Proof:} Assuming $\Rolo{af^{i}}=0$ for all $i$ we get
\bae
0 \fe \Rolo{S^{-1}(e_{i})} \Rolo{a f^{i}_{(1)}} S(f^{i}_{(2)}) \ff
 \fe \Rolo{S^{-1}(e_{i})} \Rolo{a_{(1)} f^{i}_{(1)}} a_{(2)} f^{i}_{(2)}
S(f^{i}_{(3)}) \ff
 \fe \Rolo{S^{-1}(e_{i})} \Rolo{a_{(1)} f^{i}} a_{(2)} \ff
 \fe \Rolo{S^{-1}(e_{i})} \Rolo{f^{i}} a \ff
 \fe a \,  \nn  
\eae
where, in the last line, use was made of\men{intdeltaR}.
\hfill \prend

\noindent There exist formulas similar to\men{intdefxLaR} for other
combinations of invariance properties
\bae
\la x \Lar \la a \Lar & \sim &  \ip{e_{i} x}{S^{-2}(f^{i}) a} \ff
\la x \Lar \la a \Rar & \sim & \ip{x S^{2}(e_{i})}{f^{i} a} \ff
\la x \Rar \la a \Lar & \sim & \ip{e_{i} x}{a S^{2}(f^{i})}.
\label{xLRaLR}
\eae
We could have used any of these formulas as our basic definition of
the integral. When we deal with braided Hopf algebras we use in
fact the analogue of the second of\men{xLRaLR} as our starting
point.
The coefficient of proportionality in the above formulas depends on
which function's integral we normalize to 1.
The relation with the trace formula\men{RadL} is illuminated by the
following
\begin{lemma}
\label{xtratr}
For $\A$, $\U$ dually paired FDHAs, it holds ($a$ in $\A$, $x$ in
$\U$)
\end{lemma}  
\be
\la x \trRar \la a \trRar = \sigma \ip{xe_{i}}{S^{2}(f^{i})a} 
\label{trvrel}.
\ee
{\em Proof:} 
\bae
 \la x \trRar \la a \trRar \fe \inprod{e_{i} x}{S^{2}(f^{i})}
\inprod{e_{k}}{S^{2}(f^{k}) a} \ff
 \fe \inprod{e_{i} x_{(1)}}{S^{2}(f^{i})} \inprod{x_{(2)}
e_{k}}{S^{2}(f^{k}) a} \ff
 \fe \inprod{e_{i}}{S^{2}(f^{i}_{(1)})}
\inprod{x_{(1)}}{S^{2}(f^{i}_{(2)})} \inprod{x_{(2)} e_{k}}{\Sfk )
a} \ff
 \fe \inprod{e_{i}}{\Sfi_{(1)})} \inprod{x_{(1)} {e_{k}}_{(2)}
S^{-1}({e_{k}}_{(1)})}{\Sfi_{(2)})} \inprod{x_{(2)} {e_{k}}_{(3)}}{\Sfk)
a} \ff
 \fe \inprod{e_{i}}{\Sfi_{(1)})}
\inprod{S^{-1}({e_{k}}_{(1)})}{\Sfi_{(3)})} \inprod{x
{e_{k}}_{(2)}}{\Sfi_{(2)}) \Sfk) a} \ff
 \fe \inprod{e_{i} e_{j} e_{l}}{\Sfi)} \inprod{{e_{k}}_{(1)}}{S(f^{l})} 
\inprod{x {e_{k}}_{(2)}}{\Sfj) \Sfk) a} \ff
 \fe \inprod{e_{i} e_{j} S({e_{k}}_{(1)})}{\Sfi)} \inprod{x
{e_{k}}_{(2)}}{\Sfj)) \Sfk) a} \ff
 \fe \inprod{e_{i} {e_{j}}_{(1)} S({e_{j}}_{(2)})}{\Sfi)} \inprod{x
{e_{j}}_{(3)}}{\Sfj) a} \ff
 \fe \inprod{e_{i}}{\Sfi)} \inprod{x e_{j}}{\Sfj) a} . \nn 
\eae
\hfill \prend

\noindent As a corollary, we infer that, when 
$\sigma = 0$, $\la \cdot \trRar$
 is trivial in $\U$ or $\A$ (or both); when $\sigma \neq 0$, it is
nontrivial in both $\U$ and $\A$.
\subsection{Vacuum Projectors}
\label{vacproj}
We give here a formulation of invariant integration in which the
integral of a function is regarded as its ``vacuum expectation
value''. First, notice that right invariance 
can also be expressed as ($x$ in $\U$, $a$ in $\A$)
\be
\la x \tr a \ar = \epsilon(x) \la a \ar \, . \label{rint2}
\ee
Recall now the $\U$ and $\A$-vacua introduced in section\mens{HA}. 
We could define, in terms of these, our ``vacuum'' integral 
via an equation like ($a \in \A$)~\cite{BZpc}
\be
\la a \ar_{v} \sim \ulv a \urv ;
\ee
invariance in the form\men{rint2} is automatically satisfied.
However, as we shall soon see, it is more natural, in this case, 
to work instead with
quantities like $\rva \lva$, $\rva \ulv$ \etc , \ie with operators
rather than states. The reason is that the former can be realized in
$\smash$ and hence their properties can be derived while those of the
latter have to be introduced ``by hand''. We aim therefore at a
definition like
\be
\la a \ar_{v} \sim \rva \ulv a \urv \lva \label{defofint} .
\ee
We expect the rhs of\men{defofint} to be proportional to $\da$ (at
least, under certain conditions), consistent with its
property to return counit when multiplied by
functions either from the left or from the right.
What we need next is to find quantities in $\smash$ that represent the
operators $\rva \ulv$, $\urv \lva$. We recall at this point a
result of~\cite{CSW}: the {\em vacuum projectors} $E,\, \bar{E}$ 
defined by
\be
E =  S^{-1}(f^{i}) e_{i} \, , \qquad
\bar{E} = S^{2}(e_{i}) f^{i} 
\ee
satisfy
\bae
Ea \fe S^{-1}(f^{i}) e_{i} a \ff
 \fe S^{-1}(f^{i}) a_{(1)} \ip{e_{i_{(1)}}}{a_{(2)}} e_{i_{(2)}}
\ff
 \fe S^{-1}(f^{j}) S^{-1}(f^{i}) a_{(1)} \ip{e_{i}}{a_{(2)}} e_{j}
\ff
 \fe S^{-1}(f^{j}) S^{-1}(a_{(2)}) a_{(1)} e_{j} \ff
 \fe \epsilon(a) E
\label{Ea}
\eae
for all $a$ in $\A$, as well as
\bae
xE \fe x S^{-1}(f^{i}) e_{i} \ff
 \fe S^{-1}(f^{i}_{(2)}) \ip{x_{(1)}}{S^{-1}(f^{i}_{(1)})} x_{(2)}
e_{i} \ff
 \fe S^{-1}(f^{j}) \ip{x_{(1)}}{S^{-1}(f^{i})} x_{(2)} e_{i} e_{j}
\ff
 \fe S^{-1}(f^{j}) x_{(2)} S^{-1}(x_{(1)}) e_{j} \ff
 \fe \epsilon(x) E
\eae
for all $x$ in $\U$, while we can similarly show that
\bae
\bar{E} x \fe \epsilon(x) \bar{E} \; \; \; \; \; \; \forall x \in
\U \\
a \bar{E} \fe \epsilon(a) \bar{E} \; \; \; \; \; \; \forall a \in
\A .
\eae
Furthermore, $E^{2}=E$ and $\bar{E}^{2}=\bar{E}$
 which allows us to write
$E = \urv \lva$, $\bar{E} = \rva \ulv$.
With an eye on\men{defofint}, we now compute
\bae
\bar{E} a E \fe S^{2}(e_{i}) f^{i} a E \ff
 \fe f^{i}_{(1)} a_{(1)} \ip{S^{2}(e_{i})}{f^{i}_{(2)} a_{(2)}} E
\ff
 \fe f^{n} \ip{e_{n}}{f^{i}_{(1)} a_{(1)}}
\ip{S^{2}(e_{i})}{f^{i}_{(2)} a_{(2)}} E \ff
 \fe f^{n} \ip{e_{n} S^{2}(e_{i})}{f^{i}a} E \ff
 \fe \Rolo{a} \Lda E \, . \label{EbaraE}
\eae
This simplifies further when $\Lda = \Rda \equiv \da$ - we then get
\be
\bar{E} a E = \olo{a} \da \, .
\label{EbaEda}
\ee
\subsection{Fourier Transforms}
We work again with a general (\ie not necessarily unimodular) FDHA.
We define a {\em right Fourier transform} $\ftR{\cdot} : \A \ra
\U$ in terms of a right integral as follows
\be
\ftR{a} \equiv \Rolo{aS^{-1}(f^{i})} e_{i} \label{ftdef1}
\ee
so that ($b$ in $\A$)
\be
\ip{\ftR{a}}{b} = \Rolo{aS^{-1}(b)} . \label{ftdef2}
\ee
We show now that the right Fourier transform is invertible 
\bae
\Rolo{ e_{j} \ftR{a}}  f^{j} \fe  \Rolo{a S^{-1}(f^{i})} 
\Rolo{ e_{j} e_{i} } f^{j} \ff
 \fe \Rolo{ a S^{-1}(f^{i}_{(2)}) } \Rolo{ e_{i} }
f^{i}_{(1)} \ff
 \fe \Rolo{ a_{(1)} S^{-1}(f^{i}_{(3)}) } \Rolo{ e_{i}
} a_{(2)} S^{-1}(f^{i}_{(2)}) f^{i}_{(1)} \ff
 \fe \Rolo{ a_{(1)} S^{-1}(f^{i}) } \Rolo{ e_{i} }
a_{(2)} \ff
 \fe \Rolo{ a_{(1)} S^{-1}(\Rda) } a_{(2)} \ff
 \fe \Rolo{ S^{-1}(\Rda)} a \ff
 \fe a \label{ftinv} \, .
\eae
In the language of the previous section, when $\sigma' \neq 0$,
 the Fourier transform
allows the switching between $\urv$ and $\rva$. Indeed, we find
\bae
aE \fe a S^{-1}(f^{i}) e_{i} \ff
 \fe e_{i_{(2)}} \ip{e_{i_{(1)}}}{S^{-1}(a_{(2)}
S^{-1}(f^{i}_{(1)}))} a_{(1)} S^{-1}(f^{i}_{(2)}) \ff
 \fe e_{i_{(2)}} \ip{S^{-2}(e_{i_{(1)}})}{f^{i}_{(1)} S(a_{(2)})}
a_{(1)} S^{-1}(f^{i}_{(2)}) \ff
 \fe e_{j} \ip{S^{-2}(e_{i})}{f^{i}_{(1)} f^{j}_{(1)} S(a_{(2)})}
a_{(1)} S^{-1}(f^{j}_{(2)}) S^{-1}(f^{i}_{(2)}) \ff
 \fe e_{j} \ip{S^{-2}(e_{i}) S^{-2}(e_{k})}{f^{i} f^{j}_{(1)}
S(a_{(2)})} a_{(1)} S^{-1}(f^{j}_{(2)}) S^{-1}(f^{k}) \ff
 \fe e_{j} \la S^{-2}(e_{k}) \ar \la f^{j}_{(1)} S(a_{(2)})
\ar a_{(1)} S^{-1}(f^{j}_{(2)}) S^{-1}(f^{k}) \ff
 \fe e_{j} \la a_{(2)} S^{-1}(f^{j}_{(1)}) \ar a_{(1)}
S^{-1}(f^{j}_{(2)}) \la S^{-3}(e_{k}) \ar f^{k} \ff
 \fe e_{j} \la a S^{-1}(f^{j}) \ar \la e_{k} \ar f^{k} \ff
 \fe \hat{a} \da \, , 
\eae
which, in terms of the action on right vacua, corresponds to
\be
a \urv = \hat{a} \rva \, . \label{avuftava}
\ee
One can easily check that ($a$ in $\A$, $x$ in $\U$)
\be
\ftR{x \tr a} = x \ftR{a} \, . \label{xaR}
\ee
Another familiar, in the unimodular case, property that survives
when $\olo{\uA} = 0$, is
\be
\ftR{f \star^{\scr R} g} = \ftR{f} \ftR{g} \, ,
\label{fgconvft}
\ee
where the {\em right convolution} $\star^{\scr R}$ of
$f,\, g$ in $\A$ is given by
\be
f \star^{\scr R} g = g_{(1)} \Rolo{f S^{-1}(g_{(2)})}
\label{fgconv}
\ee
(\men{fgconvft}, together with the invertibility of $\ftR{\cdot}$,
guarantee the associativity of $\star^{\scr R}$). On
the other hand, the following property that is easily seen to hold
in the unimodular case,
\be
\hat{\hat{a}} = S(a) \label{aft2Sa}
\ee
does not hold, in general, for $\ftR{\cdot}$ when $\olo{\uA} = 0$.
\subsection{Further properties}
\label{furthprop}
We give now the proof of a number of interesting formulas, valid
for unimodular FDHAs. First, notice that
\bae
\widehat{S^{2}(a)} \fe e_{i} \la S^{2}(a) S^{-1}(f^{i}) \ar \ff
 \fe e_{i} \la a S^{-3}(f^{i}) \ar \ff
 \fe S^{-2}(e_{i}) \la a S^{-1}(f^{i}) \ar \ff
 \fe S^{-2}(\hat{a}) \label{S2ft}
\eae
where, in the second line, we used\men{intaSa}. Two useful lemmas
follow 
\begin{lemma}
For $\A$ a FDHA and $\sigma' \neq 0$, it holds
\be
S^{-2}(a) =
\ip{a^{(\bar{1})}}{S^{-2}({a^{(2)}}_{(2)})} {a^{(2)}}_{(1)} 
\label{ident1}
\ee
for all $a$ in $\A$.
\end{lemma}

{\em Proof:}
We have (the notation is introduced in section\mens{HA})
\bae
a^{(\bar{1})} \ot {a^{(2)}}_{(1)} \ot {a^{(2)}}_{(2)} \fe (\id \ot
\Delta) (S(e_{i}) e_{j} \ot f^{i} a f^{j} ) \ff
 \fe S(e_{k}) S(e_{i}) e_{j} e_{l} \ot f^{i} a_{(1)} f^{j} \ot
f^{k} a_{(2)} f^{l} \nn
\eae
which gives
\bae
\inprod{a^{(\bar{1})}}{S^{-2}({a^{(2)}}_{(2)})} {a^{(2)}}_{(1)} 
\fe \inprod{S(e_{k})
S(e_{i}) e_{j} e_{l}}{S^{-2}(f^{k}) S^{-2}(a_{(2)}) S^{-2}(f^{l})}
f^{i} a_{(1)} f^{j} \ff
 \fe \la S(e_{k}) S(e_{i}) e_{j} \ar \la S^{-2}(f^{k})
S^{-2}(a_{(2)}) \ar f^{i} a_{(1)} f^{j} \ff
 \fe \la S(e_{i}) e_{j} \ar \la S^{-2}(f^{i}_{(2)}) S^{-2}(a_{(2)})
\ar f^{i}_{(1)} a_{(1)} f^{j} \ff
 \fe \la S(e_{i}) e_{j} \ar \la f^{i} a \ar f^{j} \ff
 \fe \la S(e_{i}) e_{j} \ar \la S^{-1}(a) S^{-1}(f^{i}) \ar f^{j}
\ff
 \fe \la S(\widehat{S^{-1}(a)}) e_{j} \ar f^{j} \ff
 \fe \olo{e_{j} S^{2}(\widehat{S^{-1}(a)})} S(f^{j}) \ff
 \fe \olo{e_{j} \widehat{S^{-3}(a)}} S(f^{j}) \ff
 \fe S^{-2}(a) \, ,
\eae
where, in the last line, we used the formula for the inverse
Fourier transform, eqn\men{ftinv}. \hfill \prend

\begin{lemma}
For $\A$ a FDHA, it holds
\be
a^{(2)} S^{-1}(a^{(\bar{1})}) b = b a^{(2)} S^{-1}(a^{(\bar{1})})
\label{atildeb}
\ee
for all $a, \, b$ in $\A$.
\end{lemma}
{\em Proof:} We have
\bae
 a^{(2)} S^{-1}(a^{(\bar{1})}) b \fe a^{(2)} b_{(1)}
\inprod{S^{-1}({a^{(\bar{1})}}_{(2)})}{b_{(2)}}
S^{-1}({a^{(\bar{1})}}_{(1)}) \ff
 \fe b_{(1)} \inprod{a^{(2)(\bar{1})}}{b_{(2)}} a^{(2)(2)}
\inprod{\Saot}{b_{(3)}} \Sioo \ff
 \fe b_{(1)} \inprod{a^{(2)(\bar{1})} \Saot}{b_{(2)}} a^{(2)(2)}
\Sioo \ff
 \fe b_{(1)} \inprod{{a^{(\bar{1})}}_{(3)} \Saot}{b_{(2)}} a^{(2)}
\Sioo \ff
 \fe b a^{(2)} \Sia) 
\eae
where, in the first and second line, we used the first and second 
of\men{AUcr} respectively. \hfill \prend

\noindent At this point, we have enough machinery at our disposal 
to prove the following
\begin{prop}
For $\A$ a FDHA and $\sigma' \neq 0$, it holds 
\be
 \la ba \ar = \la S^{-2}(a) b \ar
\label{bacrs}
\ee
for all $a$, $b$ in $\A$.
\end{prop}

{\em Proof:} We have
\bae
\bar{E} ba E \fe \bar{E} b a^{(2)} \epsilon (\Sia)) E \ff
 \fe \bar{E} b a^{(2)} \Sia) E \ff
 \fe \bar{E} a^{(2)} \Sia) b E \ff
 \fe \bar{E} \Sioo \inprod{\Saot}{\Stta} {a^{(2)}}_{(1)} b E \ff
 \fe \bar{E} \inprod{\Sia)}{\Stta} {a^{(2)}}_{(1)} b E \ff
 \fe \bar{E} S^{-2}(a) b E \nn
\eae
where in the third line we used\men{atildeb} and in the last
one,\men{ident1}. The proposition follows now from\men{EbaEda}.
 \hfill \prend 

\noindent It is interesting to compare\men{bacrs} with the formula
one can derive in the presence of a universal $R$-matrix $\R$ in
$\U$. The result in this case is contained in
\begin{prop}
Let $\A$ be a dual quasitriangular Hopf algebra and $\olo{\cdot} $
a biinvariant integral on it. It holds
\be
\la ba \ar = \la S^{2}(a \triangleleft s) b \ar . \label{ident3}
\ee
for all $a$, $b$ in $\A$, where 
\be
s = cu^{-2}, \, \, \, \, \, \, 
u = S(\R^{(2)}) \R^{(1)}, \, \, \, \, \, \, 
c = u S(u)
\label{sucdef}
\ee
( $s$, $u$, $c$ in $\U \sim \A^{\ast}$) and $a \triangleleft s
\equiv \ip{s_{(1)}}{a} s_{(2)}$.
\end{prop}

{\em Proof:} 
The commutation relations in $\A$ can be written in the form
\be
ba = \ip{\R}{a_{(1)} \ot b_{(1)}} a_{(2)} b_{(2)}
\ip{\R^{-1}}{a_{(3)} \ot b_{(3)}} \label{AAcrs2}
\ee
(this is the dual version of\men{RcopRi}). It can also be shown that the
element $u$  defined above implements the square of
the antipode in $\U$ acccording to
\be
S^{2}(x) = u x u^{-1} \label{uxui}
\ee
for all $x$ in $\U$~\cite{DrinfQG} - its inverse is given by $u^{-1} =
\R^{(2)} S^{2}(\R^{(1)})$. We can then write
\bae
\la ba \ar \fe \ip{\R}{a_{(1)} \ot b_{(1)}} \la a_{(2)} b_{(2)} \ar
\ip{\R^{-1}}{a_{(3)} \ot b_{(3)}} \ff
 \fe \ip{\R}{a_{(1)} \ot b_{(1)} S(b_{(2)}) 
 S(a_{(2)})} \la a_{(3)} b_{(3)} \ar
 \ip{\R^{-1}}{a_{(5)} \ot b_{(5)} S^{-1}(b_{(4)}) S^{-1}(a_{(4)})} \ff
 \fe \ip{\R}{a_{(1)} \ot S(a_{(2)})} \la a_{(3)} b 
 \ar \ip{\R^{-1}}{a_{(5)} \ot
 S^{-1}(a_{(4)})} \ff
 \fe \ip{S(u)}{a_{(1)}} \la a_{(2)} b \ar \ip{u^{-1}}{a_{(3)}} .
 \label{temp4}
\eae
However,
\bae
S(u) x u^{-1} \fe S(u)  u^{-1} u x u^{-1} \ff
 \fe s S^{2}(x)  \nn .
\eae
 Then, for arbitrary $x$ in $\U$, we have
\bae
\ip{S(u)}{a_{(1)}} \ip{x}{a_{(2)}} \ip{u^{-1}}{a_{(3)}} \fe \ip{S(u) x
u^{-1}}{a} \ff
 \fe \ip{s S^{2}(x) }{a} \ff
 \fe \ip{x}{S^{2}(a \triangleleft s)} \nn
\eae
therefore
\be
\ip{S(u)}{a_{(1)}} a_{(2)} \ip{u^{-1}}{a_{(3)}} = 
S^{2}(a \triangleleft s) ;
\label{ident2}
\ee
substituting in\men{temp4} we get\men{ident3}. \hfill \prend

\noindent Comparison with\men{bacrs} and use of\men{intnondeg} 
leads to the relation $S^{4}(a) = a \triangleleft s^{-1}$ which in 
the dual
implies $S^{4}(x) = s^{-1}x$ and therefore (by taking $x=\uU$)
\be
S^{4}(x) =x 
\label{ident4}
\ee
for all $x$ in $\U$ (a different proof of this has been given 
in~\cite{Rad2}).
Equations\men{bacrs} and\men{ident4} have been proven above only in the
finite dimensional case (the latter assuming quasitriangularity as
well). On the other hand,\men{ident3} and the
following versions of it which are proved similarly, hold for all
quasitriangular Hopf algebras with biinvariant integral
\bae
\la ba \ar \fe \la a S^{2}(s^{-1} \tr b ) \ar \ff
 \fe \la S^{-2}(s \tr a) b \ar \ff
 \fe \la (u^{-1} \tr b)(a \triangleleft u) \ar . \label{ident5}
\eae 
We close this section with the remark that $\olo{ba} = 
\olo{S^{2}(a) b}$ has been shown to hold for
unimodular, finite dimensional ribbon Hopf algebras
(see~\cite{Rad}). Using\men{ident4} in\men{bacrs} we conclude that
it actually holds for (the wider class of) quasitriangular FDHAs.
\section{Integration On Braided Hopf Algebras}
\label{IoBHA}
We transcribe here the main results of section\mens{IoHA} to the
case of FDBHAs, using the N-dimensional
quantum superplane as a concrete example.
\subsection{Preliminaries}
\subsubsection{The quantum superplane}
\label{qsuperp}
Let us review briefly the basics of the construction of the
quantum superplane~\cite{WZ,BZ2}. As is typical in the study of quantum
spaces, one deals with the associative, noncommutative algebra $\X$
generated by $1_{\X}$ and the coordinate functions 
$\xi_{i}$, $i=1,\ldots ,N$ 
on the quantum superplane satisfying the commutation relations
\be
\xi_{2} \xi_{1} = -q \hR_{12} \xi_{2} \xi_{1}    \label{xxcrs}
\ee
(we work with the ``$q^{-1}$'' version~\cite{CSZ} - we remind the reader that
$\hR_{12}(q^{-1}) = \hR_{21}^{-1}(q)$).
The derivatives $\sigmai$, $i=1, \ldots ,N$, dual to the above
coordinates, generate (together with $1_{\D}$) the algebra $\D$
with commutation relations
\be
\sigma_{1} \sigma_{2} = -q \sigma_{1} \sigma_{2} \hR_{12} .
 \label{ddcrs}
\ee
The coordinate-derivative duality is encoded in the cross relations
\be
\sigmai \xi_{j} = \delta_{ij} - q \hR^{-1}_{mj,ni} \xi_{n}
\sigma_{m} \, .
\label{dxcrs}
\ee
We denote the combined coordinate - derivative algebra by $\P$. In
analogy with the treatment of the quantum plane in~\cite{CZ}, one
can enlarge $\P$ by the introduction of displacements $\eta_{i}, \,
\tau_{i}, \, i=1, \ldots, N$ for the coordinates and derivatives
respectively, satisfying
\be
\qquad
\eta_{2} \eta_{1} = -q \hR_{12} \eta_{2} \eta_{1} \, ,
\qquad
\tau_{1} \tau_{2} = - q \tau_{1} \tau_{2} \hR_{12} \, ,
\qquad
\tau_{i} \eta_{j} = \delta_{ij} - q \hR^{-1}_{mj,ni} \eta_{n}
\tau_{m} \, ,
\qquad
\ee
\be
\qquad
\xi_{2} \eta_{1} = -q \hR^{-1}_{12} \eta_{2} \xi_{1} \, ,
\qquad
\sigma_{i} \eta_{j} = - q^{-1} \hR_{kj,li} \eta_{l} \sigma_{k} \, ,
\qquad
\label{xese}
\ee
\be
\qquad
\xi_{i} \tau_{j} = -q^{-1} D_{la} \hR_{ia,bk} D^{-1}_{bj} \tau_{l}
\xi_{k} \, ,
\qquad
\sigma_{1} \tau_{2} = -q \tau_{1} \sigma_{2} \hR^{-1}_{12} \, .
\qquad
\label{etxs}
\ee
As seen from above, the $\eta$'s are taken to be just a second copy of
the $\xi$'s but are endowed with nontrivial statistics with both the
$\xi$'s and the $\sigma$'s - analogous remarks hold for the $\tau$'s.
 The remarkable property of\men{etxs} is that the displaced
coordinates $\xi_{i} + \eta_{i}$ and derivatives $\sigma_{i} + \tau_{i}$
 still satisfy\men{xxcrs} and\men{ddcrs} respectively while the entire
enlarged algebra is covariant under the $\GL$-coaction
\bae
\xi_{i} \mapsto \xi_{i}' \fe \xi_{j} \ot S(A_{ij}) \label{xcoact} \\
\eta_{i} \mapsto \eta_{i}' \fe \eta_{j} \ot S(A_{ij}) \label{acoact} \\
\sigma_{i} \mapsto \sigma_{i}' \fe \sigma_{j} \ot S^{2}(A_{ji}) 
\label{dcoact}
\\
\tau_{i} \mapsto \tau_{i}' \fe \tau_{j} \ot S^{2}(A_{ji}) \label{pcoact} .
\eae
We will often drop the tensor product sign in the following.
\subsubsection{Braiding}
\label{braiding}
Suppose $\U$ is a quasitriangular Hopf algebra (with universal
$R$-matrix $\R$) that acts from the left on two algebras $V,\, W$.
One can, in this case, form the {\em braided tensor product} $W
\uot V$ in which $V, \, W$ are trivially embedded as subalgebras
but have nontrivial statistics, given by ($v$ in $V$, $w$ in $W$)
\bae
(1 \uot v) (w \uot 1) & \equiv & \Psi (v \ot w) \ff
 \fe \tau \circ (\R^{(1)} \tr v \ot \R^{(2)} \tr w) \ff
 \fe w^{(1)} \uot v^{(1)} \ip{\R}{v^{(2')} \ot w^{(2')}} .
\label{vwcrs}
\eae
We have expressed above the action of $\U$ on $V,\, W$ in terms of
the dual coaction of $\A \sim \U^{*}$. The
first line of\men{vwcrs} also defines the {\em braided
transposition} $\Psi: \, V \ot W \rightarrow W \ot V$ for which it
holds in general $\Psi^{2} \neq \id$ (due to $\R' \R \neq 1$). For
a detailed discussion of the properties of $\Psi$ see \eg \cite{Majid4}.
\subsubsection{The quantum superplane as a braided Hopf algebra}
\label{qpBHA}
The concept of braided tensor products provides a natural framework
for an elegant description of the quantum superplane as a braided Hopf
algebra. We give here, for completeness, an outline of this
approach - more details can be found in~\cite{Majid4}. Essential to
this description is the use of diagrams which encode neatly the
braiding information. The maps $\Psi$ and $\Psi^{-1}$ are
represented by the diagrams \, \hbox{{\epsfbox{psi.eps}}} 
\, and \, 
\hbox{{\epsfbox{psii.eps}}} \, respectively. 
For the algebra $\X$, 
the map $\xi_{i} \mapsto \xi_{i}
\uot 1 + 1 \uot \xi_{i} \equiv \eta_{i} + \xi_{i}$ is regarded as a
{\em braided coproduct } $\udel : \X \rightarrow \X \uot \X$
(extended (braided) multiplicatively on the whole $\X$).
Diagramaticaly this appears as 
  
\bigskip

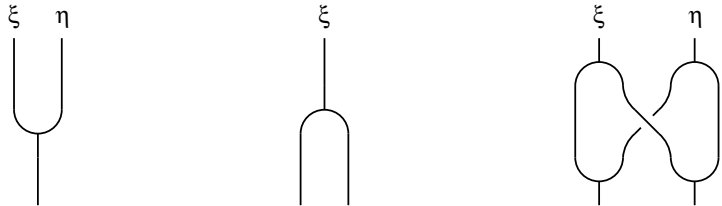
\centerline{\epsfbox{pcop.eps}}
  
\bigskip
\noindent where the first two vertices 
denote the product and coproduct in
$\X$ and the third diagram expresses the braided multiplicativity
of $\udel$ ($\xi$, $\eta$ \etc denote generic elements of $\X$).
 One also 
has a matching counit and antipode with
$\epsilon$, $S$, $S^{2}$, $S^{-1}$, $S^{-2}$ denoted respectively
by

\bigskip

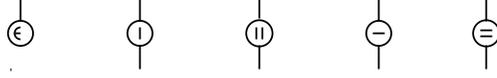
\centerline{\epsfbox{eS.eps}}
  
\bigskip

\noindent and satisfying braided versions of the familiar Hopf algebra
identities, \eg 
\be
\hbox{{\epsfbox{bpeS.eps}}}
\label{bpeS.fig}
\ee
\noindent A particularly important
requirement on the braiding, which $\Psi$ of\men{vwcrs} satisfies,
is that one should be able to move crossings past all vertices and
boxes, \eg
\be
\fgr{0}{funct.eps}
\qquad \label{funct.fig}
\ee
\noindent should hold. Exactly analogous treatment 
is possible for $\D$, the
algebra of derivatives. The braided coproduct is given by
$\udel(\sigma_{i}) = \sigma_{i} \uot 1 + 1 \uot \sigma_{i} \equiv 
\tau_{i} +
\sigma_{i}$ and the corresponding diagrams are an exact copy of those
for $\X$ (we will use in them the letters $\sigma$, $\tau$ \etc to 
denote generic elements of $\D$). One can now combine $\X$ and $\D$ 
to form a {\em braided
semidirect product} $\X \underline{\rtimes} \D$. We need for this a
braided action of $\D$ on $\X$ which is given by the second of the
following diagrams
\be
\ip{\sigma}{\xi} \, \sim \, \fgr{-1.5}{ipsixi.eps} \, , \qquad \qquad
\qquad \qquad \fgr{-1.5}{siactxi.eps}
\ee
\noindent while the first one simply 
depicts the pairing between a derivative
and a function (defined as the counit of the derivative of the
function). When viewed upside - down,
the first diagram stands for the canonical element $\phi^{i} \uot
\epsilon_{i} \in \X \uot \D$ (notice the reversal of order in the
tensor product) with
$\ip{\epsilon_{i}}{\phi^{j}}=\delta_{i}^{j}$. Both the inner
product and the canonical element are assumed invariant under
$\delA$:
\bae
\ip{p^{(1)}}{x^{(1)}} p^{(2')} x^{(2')} \fe \ip{p}{x} \uA \label{ipinv}
\\
{\phi^{i}}^{(1)} \uot {\epsilon_{i}}^{(1)} \otimes
{\phi^{i}}^{(2')} {\epsilon_{i}}^{(2')} \fe \phi^{i} \uot
\epsilon_{i} \otimes \uA \, .
\label{ceinv}
\eae
 The product - coproduct 
duality between $\X$ and $\D$ is taken to be
\be
\fgr{0}{dual.eps}
\label{dual.fig}
\ee
\noindent Notice that this differs from the standard convention in the
unbraided case. As a result, to get the unbraided version of 
any diagrammatic equation that appears in the following, one should
translate the diagrams, ignoring the braiding information, into
the language of section\mens{HA} and then set $\Delta \rightarrow
\Delta', \, S \rightarrow S^{-1}$. Again, viewing the diagrams
upside - down reveals additional (dual) information - in the case
of the diagrams above,
one discovers two basic properties of the canonical element
(compare with the first two of\men{Cprop1}). 

The commutation relations in the semidirect product (\ie the
braided analogue of\men{AUcr}) are
\be
\fgr{0}{dxcrs.eps}
\ee
We close this review with a technical remark.
 If one computes the
$\xi -\xi$ braiding given by\men{vwcrs} (with $V=W=\X$), using the
coaction\men{xcoact}, one fails to reproduce the $\xi - \eta$
commutation relations of\men{xese} - the result is off by a $q$
factor (similarly for the rest of\men{xese},\men{etxs}). To
remedy this, one can enlarge $\A$ by a grouplike central element
$g$ (the {\em dilaton}, see ~\cite{Majid2}), the inner product of
which with $\R$ is given by
\be
\ip{\R}{g^{a} \ot g^{b}} = (-q)^{-ab}, \, \, \, \, \, \, 
\ip{\R}{g^{a} \ot A_{ij}} = \ip{\R}{A_{ij} \ot g^{a}} =
\delta_{ij}.
\label{Rgip}
\ee
We will call the enlarged algebra $\tilde{\A}$. Setting $A
\rightarrow gA$ in the rhs of\men{xcoact}-\men{pcoact} gives an
$\tilde{\A}$-coaction on the quantum superplane which reproduces,
via\men{vwcrs}, the commutation relations\men{xese},\men{etxs}.
\subsection{The invariant integral}
\subsubsection{First definition and problems}
The integral we are looking for is a linear map $\langle \cdot
\rangle : \, \X \rightarrow C$ which is translationally
invariant in the following sense~\cite{WZ}
\be
\langle \sigma_{i} f(\xi) \rangle = 0, \, \, \, \, i=1, \ldots, N  
\, \, \, \, \, \, \forall f \in \X
\label{infinv} .
\ee
 An equivalent formulation of invariance is~\cite{CZ}
\be
\langle f(\xi+\eta) \rangle = 1_{\X} \langle f(\xi) \rangle 
\, \, \, \, \, \, 
\forall f \in \X
\label{fininv1}
\ee
or, in braided Hopf algebra language,
\be
f_{\underline{(1)}} \langle f_{\underline{(2)}} \rangle =
1_{\X} \langle f \rangle \, \, \, \, \, \, 
\forall f \in \X.
\label{fininv2}
\ee
Representing the integral with a rhombus, we want
it to satisfy
\be
\fgr{0}{functint.eps}
\label{functint.fig}
\ee
\noindent a requirement which, as one can easily see, cannot, in
general, be satisfied.
Indeed, we only need consider the classical fermionic line with
coordinate $\xi$ and displacement $\eta$ satisfying $\xi^{2} =
\eta^{2}=0,\, \, \xi \eta = - \eta \xi$. Taking $\xi$, $\eta$
to stand for themselves in the diagram above,
 we find that the lhs is $\eta \olo{\xi}$
(since $\olo{\xi}$, being a number, braids trivially), while the
rhs is $\eta \olo{-\xi} = - \eta \olo{\xi}$, implying $\olo{\xi}=0$
which contradicts the known Berezin result.  
We conclude that $\olo{f(\xi)}$,
assumed invariant (and nontrivial), cannot both be a number and
satisfy the property expressed in the diagram above. 
\subsubsection{An improved definition}
\label{impdef}
Our treatment, in section\mens{altform}, of the integral on FDHAs
points to a simple solution to the above problem. We recall that
there, the quantity that naturally emerged, in our algebraic
formulation in terms of the modified trace formula, was the
numerical integral $\olo{\cdot}$ {\em times} a delta function (as
in the lhs of\men{intdef3}). Motivated by this, we define a new
integral $\Lbolo{\cdot} : \, \X \rightarrow k \Rdx$ (with $k$, in
our case, the complex numbers) as follows
\be
\Lbolo{\xi} = \Lolo{\xi} \Rdx \, .
\label{bintdef1}
\ee
The output braid of the integration rhombus in our diagrams will
stand accordingly for a numerical multiple of $\Rdx$. The rhombus'
inner workings are exposed in the diagram below
\be
\fgr{0}{rhoiw.eps}
\label{rhoiw.fig}
\ee
To get the product $\Rolo{\sigma} \Lolo{\xi}$ one should pair the
output braid with $\sigma$. Whether\men{rhoiw.fig}, in its present
or a suitably modified form, applies to the infinite dimensional
case (and under what conditions), is a direction for future work
(one can easily see that, in certain such cases, the rhs
of\men{rhoiw.fig} diverges).
\begin{lemma}
The integral $\Lbolo{\cdot}$ defined by\men{rhoiw.fig} is
nontrivial for every FDBHA.
\end{lemma}
{\em Proof:}
Using $(S^{-1}
\otimes \id) \circ \Psi^{-1}(\phi^{i} \otimes \epsilon_{i})$ as
input to the braided version of $\Theta$ (shown below), we find

\fgr{0}{nontriv.eps}

\noindent from which nontriviality follows. \hfill \prend

\subsubsection{Braided vacuum projectors}
\label{bvp}
We denote in the following by $\E$, $\bE$ the
braided analogues of $E$, $\bar{E}$. They are given by
\be
\E \, = \, \fgr{-3}{vac1.eps} \, , \qquad \qquad
\bE \, = \, \fgr{-4.5}{vac2.eps}
\label{bEEb}
\ee
\noindent The proof of the (analogue of)\men{Ea} is as follows
\bigskip

\centerline{\epsfbox{eproof1.eps}}
\hfill \prend

\bigskip
\noindent As before, $\bE a \E$ is a multiple of 
$\Rdx \E$ - when $\Rdx =
\Ldx \equiv \dX$, it becomes a multiple of $\dX$. Notice that, with
$\Rdx$ defined via $\ip{\tau}{\Rdx} = \Rolo{\tau}$, one gets
$\xi \Rdx = \epsilon(\xi) \Rdx$ - the difference in the order of
multiplication, compared to\men{xdU}, is due to the second
of\men{dual.fig}.

 The explicit
computation of $\bE$ is simplified by the following identity
\be
\Psi^{-1}(\phi^{i} \otimes \epsilon_{i}) = S(u^{-1}) \tr
\epsilon_{i} \otimes \phi^{i} \, .
\label{psiice}
\ee
which is easily proved using the invariance of the canonical
element.

Inspection of our definition reveals the braided version of the
trace formula\men{RadL} for the numerical integral $\la \xi
\turRar$ -
it is given by
\be
\la \xi \turRar \, = \, \fgr{-4}{bRL.eps}
\label{bRL.fig}
\ee
The proof of invariance is as follows

\bigskip

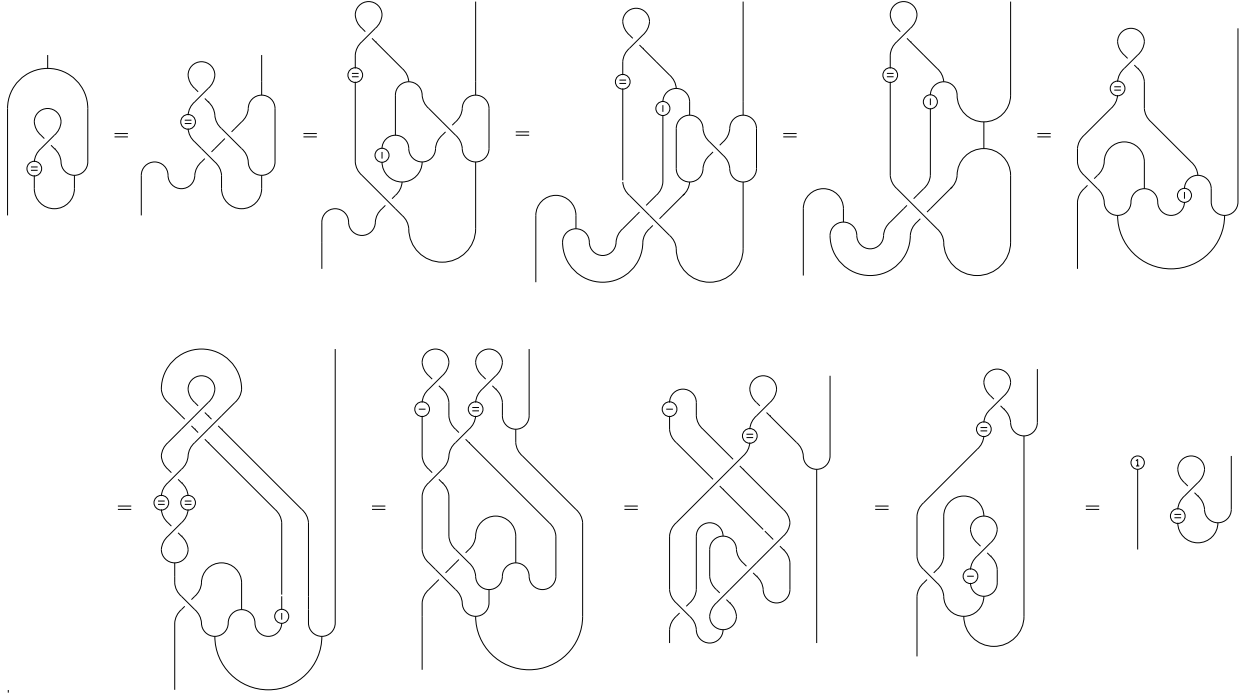
\epsfbox{invar2.eps}

\hfill \prend
  
\bigskip
To find out under what conditions it becomes trivial, we have to
derive the braided version of\men{trvrel}. Ommiting the somewhat
lengthy diagrammatic proof, which parallels that of
Lemma\mens{xtratr}, we state
\begin{lemma}
For $\X$, $\D$ dually paired FDBHAs, it holds
\end{lemma}
\be
\fgr{0}{trsixi.eps}
\label{trsixi.fig}
\ee

\bigskip

\noindent In analogy with the unbraided case, when $\underline{\sigma} 
\equiv \la 1_{\X}
\turRar = 0$, $\la \cdot \turRar$ is trivial; when
$\underline{\sigma} \neq 0$, $\la \cdot \turRar$ provides a
nontrivial integral in $\X$.
For the existence of integrals on FDBHAs and properties of them, see
also~\cite{Lyu1,Lyu2}.
An analogous definition for the numerical braided integral (and a 
different proof of its invariance) can be found in~\cite{Majid6}.
\subsubsection{Braided Fourier transforms}
Transcribing\men{ftdef1}, we define the Fourier transform $\hat{f}$ of
the element $\xi$ of a FDBHA $\X$ by the equation
\be
\hat{\xi} \equiv \olo{\xi S(\phi^{i})} \epsilon_{i}
\ee
or, in pictures,
\be
\fgr{0}{FTdef.eps}
\label{FTdef.fig}
\ee
(this differs from earlier definitions~\cite{MajKem} by the use of the
nonbosonic integral $\bolo{\cdot}$).
\noindent The output braid on the right stands for what 
one usually 
calls the Fourier transform of $\xi$ (an element of $\D$, the dual of
$\X$) while the one on the left stands for the delta function in
$\X$ that is produced by the integration and which ensures the
correct braiding behavior of $\hat{\cdot}$. There is also a notion of
braided convolution of functions, defined by
\be
\fgr{0}{convdef.eps}
\label{convdef.fig}
\ee
\noindent Again, the output braid on the left 
only carries a delta function
in $\X$. The following basic properties can be shown to hold
\bigskip

\centerline{\epsfbox{FTprop.eps}}
  
\bigskip
\subsection{Integration on the quantum superplane}
We apply now the general formalism developed
above to the problem of integration on the quantum
superplane. Our starting point will be the vacuum projector
construction of section\mens{bvp} - notice that although
$\olo{\uA} = 0$ in this case, the integral is nevertheless
biinvariant so we expect\men{EbaEda} to hold. For the canonical
element we find
\be
\phi^{i} \uot \epsilon_{i} = e_{q^{-1}}(\xi_{i} \uot \sigma_{i}) 
\,  \label{ceqsp1}
\ee
where 
\be
e_{q}(x) = \sum_{k=0}^{\infty}\frac{1}{[k]_{q}!} x^{k} \, , \qquad
[k]_{q}=\frac{1-q^{2k}}{1-q^{2}} \, , \qquad [k]_{q}! = [1]_{q}
[2]_{q} \ldots [k]_{q} \, , \qquad
[0]_{q} \equiv 1 \label{eqxdef}
\ee
(compare with the vacuum projectors for the
quantum plane in~\cite{CZ}). The commutation
relations\men{xxcrs},\men{ddcrs} imply $\xi_{i}^{2} =
\sigma_{i}^{2} =0$ for $i=1, \ldots ,N$ which gives $(\xi_{i} \uot
\sigma_{i})^{N+1} =0$. Using the braiding relations\men{xese}, the
second of which can also be written as
\be
\Psi(\sigma_{i} \ot \xi_{j}) = -q^{-1} \hR_{kj,li} \xi_{l} \ot
\sigma_{k} \, , \label{xese2}
\ee
we can expand $(\xi_{i} \uot \sigma_{i})^{k}$ in\men{ceqsp1} to
find
\be
\phi^{i} \uot \epsilon_{i} = \sum_{k=0}^{N}
\frac{q^{-k(k-1)}}{[k]_{q^{-1}}!} \xisi \, .
\label{ceqsp2}
\ee
With the antipode being given by 
\be
S(\xs) = (-1)^{k} q^{k(k-1)} \xs
\label{Ssx}
\ee
and using\men{psiice} we find
\bae
\E \fe \sum_{k=0}^{N} \frac{(-1)^{k}}{[k]_{q^{-1}}!} \xs \sx \\
\bE \fe \sum_{k=0}^{N} \frac{(-1)^{k} q^{k}}{[k]_{q}!} D_{i_{1}j_{1}}
\ldots D_{i_{k}j_{k}} \sx \xi_{j_{1}} \ldots \xi_{j_{k}} \, .
\label{bEEbqsp}
\eae
It will be convenient in the following to express $\bE$ in the
alternative form 
\be
\bE = \sum_{k=0}^{N} \frac{(-1)^{k} q^{k(k-2N+1)}}{[k]_{q}!}
([N]_{q} - \xi \! \cdot \! \sigma)([N-1]_{q} - \xi \! \cdot \! \sigma) \ldots
([N-k+1]_{q} - \xi \! \cdot \! \sigma) \label{be2}
\ee
where $\xi \! \cdot \! \sigma \equiv \xi_{i} \sigma_{i}$ (this form makes
the invariance under $\delA$ evident - a similar expression exists
for $\E$). Using the
commutation relation $\xi \! \cdot \! \sigma \xi_{j} = \xi_{j} (1 + q^{2}
\xi \! \cdot \! \sigma)$ and the fact that $\xi \! \cdot \! \sigma \E =0$, we
can now compute the integral of an arbitrary monomial $\xi_{i_{1}}
\ldots \xi_{i_{r}}$ ($r < N$)
\be
\bE \xi_{i_{1}} \ldots \xi_{i_{r}} \E =
(\sum_{k=0}^{A}\frac{(-1)^{k} q^{k(k-2A+1)} [A]_{q}!}{[k]_{q}!
[A-k]_{q}!}) \xi_{i_{1}} \ldots \xi_{i_{r}} \E
\ee
where $A \equiv N-r$. For the sum in parentheses, one can set
\be
S(z) = \sum_{k=0}^{A}\frac{(-1)^{k} q^{k(k-2A+1)} [A]_{q}!}{[k]_{q}!
[A-k]_{q}!} z^{k} \, .
\label{Szdef}
\ee
Introducing a Jackson derivative $\partial_{z}$, satisfying
$\partial_{z} z = 1 + q^{2} z \partial_{z}$, we find
from\men{Szdef}
\be
\partial_{z} S(z) = \frac{1}{q^{-2}-1} S(z) -
\frac{q^{-2A}}{q^{-2}-1} S(q^{2}z) \, .
\label{dzSz1}
\ee
On the other hand, it holds
\[
\partial_{z} S(z) = \frac{S(q^{2}z) -S(z)}{q^{2}-1} \, ;
\]
comparison with\men{dzSz1} shows that
\[
\frac{S(q^{2}z)}{S(z)} = \frac{1-q^{2}z}{1-q^{-2(A-1)}z}
\]
from which we find 
\be
S(z) = (1-z)(1-q^{-2}z) \ldots (1-q^{-2(A-1)}z) \, ,
\label{Szfinal}
\ee
implying that $S(1)=0$ and therefore that
\be
\olo{\xi_{i_{1}} \ldots \xi_{i_{r}}} =0 \qquad 0 \leq r < N \, .
\label{xiint}
\ee
For $r=N$, the integrand is (a multiple of) a delta function and 
its numerical 
integral is evidently (a multiple of) 1. We conclude that the 
quantum Berezin
integral in $N$ dimensions is essentially undeformed. As expected, 
$\la \cdot
\turRar$ is trivial in this case, as the reader can easily verify.
\subsection{Remarks on the integral on the quantum plane}
We make here a few remarks about the transformation properties of
the invariant integral on the quantum plane~\cite{WZ}. We denote by
$x_{i}, \, i=1, \ldots ,N$ the coordinate functions on it and by
$\partial_{i}, \, i=1, \ldots ,N$ the dual derivatives, satisfying
\be
x_{1} x_{2} = q^{-1} \hR_{12} x_{1} x_{2} \, ,
\qquad
\partial_{2} \partial_{1} = q^{-1} \hR_{12} \partial_{2}
\partial_{1}\, ,
\qquad
\partial_{i} x_{j} = \delta_{ij} + q \hR_{jl.ik}
x_{k} \partial_{l} \, .
\label{qpcrs}
\ee
The above algebra is covariant under the transformation $x
\rightarrow \delA(x) =xA, \, \p \ra \delA(\p) = \p S(A^{T})$ with
$A$ a $\GL$ matrix (we ommit the tensor product symbol). We
assume an integral $\bolo{\cdot}$ exists, defined on a suitable
class of functions, satisfying translational invariance (in the
spirit of\men{infinv}) and braiding correctly, \ie according
to\men{functint.fig}.
In the classical case of integration on the N - dimensional plane,
one finds the transformation property
\be
\int f(xA) dx = \frac{1}{\mathop{\rm det} (A)} \int f(x) dx \nn
\ee
where $A$ is a $GL(N)$ matrix. We now show that a similar property
holds in the quantum case. We remark first that\men{functint.fig}
implies
\be
\delA(\bolo{f(x)}) = \bolo{f(xA)} = \bolo{f^{(1)}(x)} f^{(2')}(A)
\label{intcoact}
\ee
where $f(xA) \equiv f^{(1)}(x) f^{(2')}(A)$. Consider now the dual
action of the generators $Y_{ij} \equiv L^{+}_{im} S(L^{-}_{mj})$ 
of $\U_{q}(gl(N))$ on the integrand
\bae
\bolo{ Y_{ij} \tr f} \fe \bolo{ f^{(1)}}  \ip{Y_{ij}}{f^{(2')}} .
\label{Yactf}
\eae
The above action can be represented in terms of
differential operators on the plane as follows~\cite{CSZ,ChuZ}
\be
Y_{ij} \sim q^{-2} \delta_{ij} + q^{-1} \lambda \p_{i} x_{j} 
\label{Yrepdx}
\ee
which gives, making use of the invariance of the integral
\bae
\bolo{Y_{ij} \tr f } \fe \bolo{(q^{-2} \delta_{ij} + 
q^{-1} \lambda \p_{i} x_{j}) f}
\ff
 \fe q^{-2} \delta_{ij} \bolo{f} . \nn
\eae
Repeating the calculation for products of $Y$'s acting on $f$, we
find
\be
\bolo{f^{(1)}} \ot f^{(2')} = \bolo{f} \ot z 
\ee
with $\Delta(z) = z \ot z$ and $\ip{Y_{ij}}{z} = q^{-2}
\delta_{ij}$, $\ip{\uU}{z} = 1$. The above information completely
determines $z$. Since $\ip{Y_{ij}}{\detq(A)} = q^{2}$,
$\epsilon(\detq(A))=1$ and $\detq(A)$ is grouplike, we conclude
\be
\bolo{f(xA)} = \bolo{f(x)} (\detq(A))^{-1} .
\label{Acoactint}
\ee
As in the case of the quantum superplane, to obtain the correct $q$
factors for the braiding (so that, for example, 
$x \ra x \uot 1 + 1 \uot x$ is a homomorphism of the first
of\men{qpcrs}) one has to introduce a grouplike, central dilaton $g$
(extending this way $\A$ to $\tilde{\A}$) with
\be
\ip{\R}{g^{a} \ot g^{b}} = q^{ab} \, , 
\qquad
\ip{\R}{g^{a} \ot A_{ij}} = \ip{\R}{A_{ij} \ot g^{a}} = \delta_{ij}
\label{dilip}
\ee
and use the $\tilde{\A}$ coaction $x \ra xAg$, $\p \ra \p S(A^{T})
g^{-1}$ in\men{vwcrs}. We obtain in this way the braiding relations
\bae
\Psi(\bolo{f(x)} \ot x_{i}) \fe q^{-(N+1)} x_{i} \ot \bolo{f(x)}
\ff
\Psi (\bolo{f(x)} \ot \p_{i} ) \fe q^{N+1} \p_{i} \ot \bolo{f(x)}
\ff
\Psi(\bolo{f(x)} \ot \bolo{g(x)}) \fe q^{N(N+1)} \bolo{g(x)} \ot
\bolo{f(x)} \, .
\label{bintqp}
\eae
We point out that a translationally invariant integral on the
quantum plane cannot be also invariant under the coacting quantum
group - an assumption to the contrary is made in~\cite{MajKem}.
\section{Integration on Quantum Group Modules}
We present here an approach to integration on quantum spaces that
are covariant under a quantum group transformation, which does not
rely on the braided Hopf algebra structure we have assumed so far.
The only necessary ingredients are
\begin{itemize}
\item a coaction $\delA : \X \rightarrow \X \ot \A$, where $\A$,
$\X$ are the algebras of functions on the quantum group and the
quantum space respectively 
\item a map $\eta: \X \rightarrow C$ that respects the algebra
structure of $\X$
\item a (left) invariant integral on $\A$.
\end{itemize}
Given the above data, an invariant (under $\delA$) integral on $\X$
can be defined by
\be
\olo{\alpha} = \eta(\alpha^{(1)}) \olo{\alpha^{(2')}}
\label{invintcomod}
\ee
(we denote by $\olo{\cdot}$ the integral on both $\X$ and $\A$).
Notice that our notion of invariance in this section is different
from the one employed so far in this paper. Indeed, in the absence
of a (possibly braided) Hopf structure, no concept of translation
exists (as codified by the coproduct) and therefore\men{rint}
cannot serve as our starting point. 
We illustrate the above procedure taking for $\X$ the quantum Euclidean
space (for detailed treatments of this case
see~\cite{Fiore,Stein}). An example of a function algebra the integral
 on which eludes all methods presented so far in this paper
 appears in~\cite{CHZ}).
\subsection{Integration on the quantum Euclidean space}
In the following we use the notation and conventions of~\cite{FRT}.
 The algebra of functions on the $N$-dimensional quantum Euclidean 
space is generated by the coordinates $x_{i}, \, \, i=1, \ldots , N$
satisfying 
\be
x_{1} x_{2} P^{(-)}_{12} = 0,
\label{qesxxcrs}
\ee
where $P^{(-)}$ is the antisymmetric projector in the spectral
decomposition of the $SO_{q}(N)$ $\hR$ matrix. The
center of the algebra is generated by 1 and the squared length
$L^{2} = x^{T}Cx$ where $C$ is the quantum metric. For $A$ an
$SO_{q}(N)$ matrix, it holds
\be
S^{2}(A) = DAD^{-1}, \qquad \qquad 
D=CC^{T}, \qquad \qquad
A= C^{T} S(A^{T}) (C^{-1})^{T} \, .
\label{S2ASON}
\ee
The algebra of the $x$'s admits the coaction $\delA: \, x 
\mapsto xA$ while
the map $\eta: x_{i} \mapsto u_{i} \equiv u_{1} \delta_{i1} + u_{N}
\delta_{iN}$, with $u_{1}$, $u_{N}$ numbers,
respects\men{qesxxcrs}. $L^{2}$ is $\delA$-invariant.

The integral $\I_{i_{1} \ldots i_{m}} \equiv \olo{x_{i_{1}} \ldots
x_{i_{m}}}$ is given by
\be
\I_{i_{1} \ldots i_{m}} = \eta(x_{j_{1}} \ldots x_{j_{m}})
\olo{A_{j_{1} i_{1}} \ldots A_{j_{m} i_{m}}}
\label{Iqesdef}.
\ee
Notice that the function $x_{j_{1}} \ldots x_{j_{m}} \olo{A_{j_{1}
i_{1}} \ldots A_{j_{m} i_{m}}}$ is $\delA$-invariant. Assuming that
all $\delA$-invariant functions are functions of the invariant
length, we conclude
\be
x_{j_{1}} \ldots x_{j_{m}} \olo{A_{j_{1} i_{1}} \ldots A_{j_{m}
i_{m}}} = \left \{ 
\begin{array}{cl}
f_{i_{1} \ldots i_{m}}(L^{m}) & \mbox{$m$ even} \\
0 & \mbox{$m$ odd} 
\end{array} \right. 
\ee
and therefore
\be
\I_{i_{1} \ldots i_{m}} = \left \{ \begin{array}{cl}
\eta(f_{i_{1} \ldots i_{m}}(L^{2})) & \mbox{$m$ even} \\
0 & \mbox{$m$ odd} \end{array} \right. .  
\ee
We treat, as an example, the case $m=2$. Setting
\be
F_{ni,mj} \equiv \olo{S(A_{im}) A_{nj}}
\ee
and invoking invariance in the form $S(a_{(1)}) \olo{a_{(2)}} = 1
\olo{a}$, we get
\bae
A_{rb} F_{bi,aj} \fe A_{rb} S(A_{bm}) S^{2}(A_{ni}) F_{ma,nj} \ff
\Rightarrow A_{1} F_{12} D_{1} \fe F_{12} D_{1} A_{1} .
\eae
Since only multiples of the identity matrix commute with $A$, we
conclude
\be
F_{12} = M_{2} D_{1}^{-1} \nn .
\ee
However, 
\bae
A_{ki} F_{bi,aj} \fe A_{ki} S(A_{in}) A_{rj} F_{bn,ar} \ff
\Rightarrow A_{2} F_{12} \fe F_{12} A_{2} \ff
\Rightarrow M_{2} \fe \rho I_{2} \, ,
\eae
with $\rho$ a number. With the normalization $\olo{\uA} = 1$, we
find $\rho=(\mathop{\rm Tr} D^{-1})^{-1}$ and therefore
\be
F_{bi,aj} =\frac{1}{\mathop{\rm Tr} D^{-1}} D^{-1}_{ba} \delta_{ij}
\label{Fbiaj} .
\ee
Using now the third of\men{S2ASON}
in\men{Iqesdef} and substituting the above result, we find, with
the help of\men{S2ASON},
\be
\I_{i_{1} i_{2}} = \frac{\eta (L^{2})}{\mathop{\rm Tr} D^{-1}}
C_{i_{1} i_{2}} .
\label{Ii1i2}
\ee

\subsection*{Acknowledgements}
Many people have helped me understand the subject matter treated
here. I thank Joris Van der Jeught for discussions and computations
that convinced me to treat the infinite dimensional case elsewhere.
While this paper was being written up, A. Van Daele kindly pointed
out reference~\cite{VD1} which contains
essentially\men{intdefxLaR}. Also, the alternative proof of Lemma 1
quoted here (ending in eqn.\men{intdeltaR}) motivated the proof of
Lemma 6 - I thank him for his comments during the discussions we
had in Warsaw. I also thank for discussions Shahn Majid and
Volodimir Lyubashenko (who also drew my attention to
references~\cite{Lyu1,Lyu2}). Thanks go to Daniel Arnaudon for his
help with all things binary. Special thanks go to Stanislaw
Zakrzewski for his interest in this paper and for the 
organization of (and
invitation to) the minisemester on Quantum Groups and Quantum
Spaces (November - December 1995, Warsaw) during which many aspects
of the subject treated here were illuminated. Bruno Zumino had a
decisive impact on this work through numerous discussions,
contributions and his constant support - my warmest thanks go to him. 


\end{document}